\documentstyle[12pt]{article}
\def\Bbb R{{\rm \bf R}}
\def\proclaim#1{\vskip2mm{\bf #1}\em}
\def\endproclaim{\em \vskip2mm}
\def\tag#1{\eqno(#1)}
\def\gathered{\begin{array}{c}}
\def\endgathered{\end{array}}
\def\text{\mbox}
\begin{document}

\title {The framework of the enclosure method with dynamical data
and its applications}
\author{Masaru IKEHATA\footnote{
Department of Mathematics,
Graduate School of Engineering,
Gunma University, Kiryu 376-8515, JAPAN}}
\date{   }
\maketitle
\begin{abstract}
The aim of this paper is to establish the {\it framework} of the enclosure method for some class of
inverse problems whose governing
equations are given by parabolic equations with discontinuous coefficients.

The framework is given by considering a concrete inverse initial boundary value problem for
a parabolic equation with discontinuous coefficients.
The problem is to extract information about the location and shape of unknown inclusions
embedded in a known isotropic heat conductive body from a set of the input heat flux across the boundary
of the body and output temperature on the same boundary.
In the framework the original inverse problem is
reduced to an inverse problem whose governing equation have a {\it large parameter}.
A list of requirements which enables one to apply the enclosure method to the
reduced inverse problem is given.

Two new results which can be considered as the applications of the framework are given.
In the first result the background conductive body is assumed to be {\it homogeneous}
and a family of explicit {\it complex} exponential solutions are employed.
Second an application of the framework to inclusions in an isotropic {\it inhomogeneous} heat conductive body
is given.  The main problem is the construction of the special solution of the governing equation
with a large parameter for the background inhomogeneous body required by the framework.
It is shown that, introducing another parameter which is called the {\it virtual slowness} and making it {\it sufficiently large},
one can construct the required solution which yields an extraction formula of the {\it convex hull}
of unknown inclusions in a known isotropic inhomogeneous conductive body.

\noindent
AMS: 35R30


\end{abstract}


\section{Introduction}

The aim of this paper is to establish the {\it framework} of the enclosure method \cite{IE}
for possible application to some class of inverse problems whose governing
equations are given by parabolic equations
with discontinuous coefficients.

The framework is given by considering a concrete inverse initial boundary value problm
for a parabolic equation with discontinuous coefficients.
The problem is to extract information about the location and shape of unknown inclusions
embedded in a known isotropic heat conductive body from a set of the input heat flux across the boundary
of the body and output temperature on the same boundary.

Let $\Omega$ be a bounded domain of $\Bbb R^n$, $n=~2,3$ with a smooth boundary.
We denote the unit outward normal vectors to $\partial\Omega$ by the symbol $\nu$.
Let $T$ be an arbitrary {\it fixed} positive number.

Given $f=f(x,t),\,(x,t)\in\partial\Omega\times\,]0,\,T[$
let $u=u_f(x,t)$ be the solution of the initial boundary
value problem for the parabolic equation:
$$
\begin{array}{c}
\displaystyle\partial_tu-\nabla\cdot\gamma\nabla u=0\,\,\text{in}\,\Omega\times\,]0,\,T[,\\
\\
\displaystyle
\gamma\nabla u\cdot\nu=f\,\,\text{on}\,\partial \Omega\times\,]0,\,T[,\\
\\
\displaystyle
u(x,0)=0\,\,\text{in}\,\Omega,
\end{array}
\tag {1.1}
$$
where $\gamma=\gamma(x)=(\gamma_{ij}(x))$ satisfies

$\quad$

\noindent (G1)  for each $i,j=1,\cdots,n$ $\gamma_{ij}(x)$ is real
and satisfies $\gamma_{ij}(x)=\gamma_{ij}(x)\in
L^{\infty}(\Omega)$;

$\quad$

\noindent
(G2)  there exists a positive constant $C$ such that $\gamma(x)\xi\cdot\xi\ge C\vert\xi\vert^2$ for all $\xi\in\Bbb R^n$
and a. e. $x\in\Omega$.

$\quad$

See \cite{DL} for the notion of the weak solution.
This paper is concerned with the extraction of information about
``discontinuity'' of $\gamma$ from $u$ and $\gamma\nabla
u\cdot\nu$ on $\partial\Omega\times ]0,\,T[$ for some $f$ and an
arbitrary fixed $T<\infty$. However, we do not consider completely
general $\gamma$.
Instead we assume that there exists an open set
$D$ with a smooth boundary such that $\overline D\subset\Omega$
and $\gamma(x)$ a.e. $x\in\Omega\setminus D$ coincides with
a smooth positive function $\gamma_0(x)$ of $x\in\overline\Omega$
and satisfies one of the following two conditions:

$\quad$

\noindent
(A1)  there exists a positive constant $C'$ such that $-(\gamma(x)-\gamma_0(x)I_n)\xi\cdot\xi\ge C'\vert\xi\vert^2$
for all $\xi\in\Bbb R^n$ and a.e. $x\in D$;

$\quad$

\noindent
(A2)  there exists a positive constant $C'$ such that $(\gamma(x)-\gamma_0(x)I_n)\xi\cdot\xi\ge C'\vert\xi\vert^2$
for all $\xi\in\Bbb R^n$ and a.e. $x\in D$.

$\quad$

Write $h(x)=\gamma(x)-\gamma_0(x)I_n$ a.e. $x\in D$.

We consider

$\quad$

{\bf\noindent Inverse Problem 1.1.}  Fix a $T>0$.
Assume that both $D$ and $h$ are {\it unknown}
and that $\gamma_0$ is {\it known}.
Extract information about the location and shape
of $D$ from a set of the pair of temperature $u_f(x,t)$ and
heat flux $f(x,t)$ for $(x,t)\in\,\partial\Omega\,\times\,]0,\,T[$.

$\quad$

The $D$ is a model of the union of unknown {\it inclusions} where
the heat conductivity is anisotropic, different from that of the surrounding {\it inhomogeneous} isotropic conductive medium.
The problem is a mathematical formulation of a typical inverse problem in thermal
imaging. Note that in \cite{EI} a uniqueness theorem with {\it infinitely many}
$f$ for Inverse Problem 1.1 has been established
provided $h$ has the form $bI_n$ with a smooth function $b$ on $\overline D$.
Thus the point is to give a concrete procedure or formula which yield
information about the location and shape of $D$.
Note that when $n=1$, there are some results: the procedure in \cite{DKN, DLLN}
with {\it infinitely many} $f$ and the formula of Theorem 2.1 in \cite{ION} with a {\it single} $f$.

In \cite{IKII} we considered Inverse Problem 1.1 in the case when $\gamma_0\equiv 1$ and $n=3$
and gave {\it four} extraction formulae of some information including the convex hull of $D$.
In Subsection 1.3 we will reconsider three results of those which employ infinitely many
$f$, from the view point of the framework given here.

Note that in this paper we do not consider the {\it single input} case.
Formulae with a single $f$ see Theorem 1.1 in \cite{IKII} and
\cite{IKIII} for a cavity which is the extremal case $\gamma\equiv
0$ in $D$.  Those can be considered as some kind of extension of the enclosure method
for elliptic equations with a single Cauchy data started in \cite{IS} to the parabolic equations.

\subsection{A reduction to an inverse boundary value problem with a parameter-A general framework}

Define
$$\displaystyle
w_f(x,\tau)=\int_0^T e^{-\tau t}u_f(x,t)dt,\,\,x\in\Omega,\,\tau>0.
$$
The $w=w_f$ satisfies
$$\begin{array}{c}
\displaystyle
(\nabla\cdot\gamma\nabla-\tau)w
=e^{-\tau T}u_f(x,T)\,\,\text{in}\,\Omega,\\
\\
\displaystyle
\gamma\nabla w\cdot\nu
=\int_0^T e^{-\tau t}f(x,t)dt\,\,\text{for}\,x\in\partial\Omega.
\end{array}
\tag {1.2}
$$

This motivates a formulation of the reduced problem given below.

Given $F(x,\tau)$ and $g(x,\tau)$ with let $w=w(x,\tau)$ be the solution of
$$\begin{array}{c}
\displaystyle
(\nabla\cdot\gamma\nabla-\tau)w=e^{-\tau T}F(x,\tau)\,\,\text{in}\,\Omega,\\
\\
\displaystyle
\gamma\nabla w\cdot\nu=g(x,\tau)\,\,\text{on}\,\partial\Omega.
\end{array}
\tag {1.3}
$$

{\bf\noindent Inverse problem 1.2.}
Assume that $F(x,\tau)$, $D$ and $h$ are all unknown and that $\gamma_0$ is known.
Extract information about the location and shape
of $D$ from a set of the pair of $w(x,\tau)$ and
$g(x,\tau)$ for $x\in\,\partial\Omega$, $\tau>0$.

For this problem we propose the following general framework which reduces the problem to
construct a family of special solutions of an equation coming from the background body.

\proclaim{\noindent Theorem 1.1.}
Assume that: there exist constants $C_1$ and $\mu_1$ such that,
as $\tau\longrightarrow\infty$
$$\displaystyle
\Vert F(\,\cdot\,,\tau)\Vert_{L^2(\Omega)}
=O(e^{C_1\tau}\tau^{\mu_1});
\tag {1.4}
$$
we have a family $(v_{\tau})$ indexed with $\tau\ge\tau_0>0$ of
solutions of the equation
$$\displaystyle
\nabla\cdot\gamma_0\nabla v-\tau v=0\,\,\text{in}\,\,\Omega
\tag {1.5}
$$
satisfying the conditions, for some constants $\mu_2$, $\mu_3$, $\mu_4$, $C_2$ and $C_3$
$$\displaystyle
\Vert\nabla v_{\tau}\Vert_{L^2(D)}=O(e^{C_2\tau}\tau^{\mu_2}),
\tag {1.6}
$$
$$\displaystyle
\Vert\nabla v_{\tau}\Vert_{L^2(D)}\ge C'' e^{C_2\tau}\tau^{\mu_3},
\tag {1.7}
$$
$$
\displaystyle
\Vert v_{\tau}\Vert_{H^1(\Omega)}
=O(e^{C_3\tau}\tau^{\mu_4}).
\tag {1.8}
$$

Let $g=g(x,\tau)$ be a function of $x\in\partial\Omega$ having the form
$$\displaystyle
g=\Psi(\tau)\,\gamma_0\frac{\partial v_{\tau}}{\partial\nu}\vert_{\partial\Omega},
$$
where $\Psi$ satisfies the conditions, for constants $\mu$ and $\mu'$
$$\displaystyle
\liminf_{\tau\longrightarrow\infty}\tau^{\mu}\vert\Psi(\tau)\vert>0
\tag {1.9}
$$
and
$$\displaystyle
\vert\Psi(\tau)\vert=O(\tau^{\mu'}).
\tag {1.10}
$$

Let $w$ be the solution of (1.3).  If $T$ satisfies
$$\displaystyle
T>C_1+C_3-2C_2,
\tag {1.11}
$$
then
$$\displaystyle
\lim_{\tau\longrightarrow\infty}
\frac{1}{2\tau}
\log
\left\vert\int_{\partial\Omega}
\left(g\overline v_{\tau}-w\gamma_0\frac{\partial\overline v_{\tau}}{\partial\nu}\right)dS\right\vert
=C_2.
\tag {1.12}
$$

\endproclaim

{\bf\noindent Remark 1.1.}
It follows from (1.6) and (1.7) that
$$\displaystyle
\lim_{\tau\longrightarrow\infty}\frac{1}{\tau}\log\Vert\nabla v_{\tau}\Vert_{L^2(D)}=C_2.
$$
This is the meaning of $C_2$ which is uniquely determined by
$\Vert\nabla v_{\tau}\Vert_{L^2(D)}$ with all $\tau>>\tau_0$.

{\bf\noindent Remark 1.2.}  It follows from (1.7) and (1.8) that $C_3\ge C_2$.

{\bf\noindent Remark 1.3.}  From the proof one obtains also the order of the convergence
of the formula (1.12):
$$\displaystyle
\frac{1}{2\tau}
\log
\left\vert\int_{\partial\Omega}
\left(g\overline v_{\tau}-w\gamma_0\frac{\partial\overline v_{\tau}}{\partial\nu}\right)dS\right\vert
=C_2+O\left(\frac{\log\tau}{\tau}\right).
$$
This is important for a suitable choice of $\tau$ in the case when the data is noisy.

Here we present an application of Theorem 1.1 to the case when
$\gamma_0\equiv 1$. Let $c>0$ and $\tau\ge \tau_0=c^{-2}$. Let
$\omega,\omega^{\perp}\in S^{n-1}$, $n\ge 2$ and satisfy
$\omega\cdot\omega^{\perp}=0$. Set
$$\displaystyle
z=c\tau\left(\omega+i\,\sqrt{1-\frac{1}{c^2\tau}}\,\omega^{\perp}\right).
\tag {1.13}
$$
$z$ satisfies
$$\displaystyle
z\cdot z=\tau.
\tag {1.14}
$$
We observe that
$$\displaystyle
(\triangle-\tau)e^{x\cdot z}=0.
\tag {1.15}
$$
Thus we take the family $(v_{\tau})_{\tau\ge\tau_0}$ in Theorem 1.1
$$\begin{array}{c}
\displaystyle
v_{\tau}(x)=v(x;z)=e^{x\cdot z}.
\end{array}
$$
Note that:

$\bullet$  it follows from (1.15) that the function $e^{-\tau t}v(x;z)$ of $(x,t)$ satisfies the {\it backward heat equation}
$$\displaystyle
(\triangle+\partial_t)(e^{-\tau t}v(x;z))=0.
$$

$\bullet$  the absolute value of $e^{-\tau t}v(x;z)$ coincides
with $e^{-\tau(t-cx\cdot\omega)}$ and this is a solution of the
{\it wave equation} with the propagation speed $1/c$.  By this reason we
call this $c$ the {\it virtual slowness}.  See also \cite{ION} for
this interpretation.

Recall the support functions of $D$ and $\Omega$:
$$\displaystyle
h_D(\omega)=\sup_{x\in D}x\cdot\omega,\,\,h_{\Omega}(\omega)=\sup_{x\in\Omega}x\cdot\omega.
$$
We have $h_D(\omega)<h_{\Omega}(\omega)$ for all $\omega$ since $\overline D\subset\Omega$.

Applying Theorem 1.1 to (1.2), we obtain the following corollary.

\proclaim{\noindent Corollary 1.1.}
Assume that $\gamma_0\equiv 1$.
Let $f$ be the function of $(x,t)\in\partial\Omega\times]0,\,T[$ having a parameter $\tau>0$
defined by the equation
$$\displaystyle
f(x,t)=\frac{\partial v_{\tau}}{\partial\nu}(x)\varphi(t),
\tag {1.16}
$$
where a real-valued function
$\varphi\in L^2(0,\,T)$
satisfying the condition:
there exists $\mu\in\Bbb R$ such
that
$$\displaystyle
\liminf_{\tau\longrightarrow\infty}
\tau^{\mu}
\left\vert\int_0^Te^{-\tau t}\varphi(t)dt\right\vert>0.
\tag {1.17}
$$
Let $u_f=u_f(x,t)$ be the weak solution of (1.1) for $f=f(x,t)$.
If $T$ satisfies
$$\displaystyle
T>2c(h_{\Omega}(\omega)-h_D(\omega)),
\tag {1.18}
$$
then
$$\displaystyle
\lim_{\tau\longrightarrow\infty}
\frac{1}{2\tau}
\log\left\vert
\int_{\partial\Omega}
\int_0^Te^{-\tau t}\left(-\overline{v_{\tau}(x)}f(x,t;\tau)+u_f(x,t)\frac{\partial\overline{v_{\tau}}}{\partial\nu}(x)\right)dtdS\right\vert
=h_D(\omega).
\tag {1.19}
$$

\endproclaim

{\bf\noindent Remark 1.4.}
There is no restriction on the position of center of coordinates relative to $\Omega$ and $D$
which affects on the sign of $h_D(\omega)$.

Since $-h_{\Omega}(-\omega)=\inf_{x\in\Omega}x\cdot\omega<\inf_{x\in D}x\cdot\omega\le h_D(\Omega)$, we have
$$\displaystyle
h_{\Omega}(\omega)-h_D(\omega)<h_{\Omega}(\omega)+h_{\Omega}(-\omega).
$$
Thus $T$ in Corollary 1.1 can be {\it arbitrary small} by choosing a small known $c$ in such a way that, for example,
$$\displaystyle
c\le\frac{T}{2(h_{\Omega}(\omega)+h_{\Omega}(-\omega))}
$$
since from this one gets (1.18).
Therefore if one wants to estimate $D$ from the direction $\omega$ by the formula (1.19)
and the `size' of $\Omega$ at the $\omega$ direction, that is the quantity $h_{\Omega}(\omega)+h_{\Omega}(-\omega)$,
is too large compared with $T$ (this is the most difficult case), the virtual slowness $c$ in (1.13) should be chosen very small.
This is one of the two roles of virtual slowness. In the next subsection we give another role of virtual slowness.

Let us explain how to deduce (1.19) from Theorem 1.1.
Comparing (1.2) with (1.3), one knows that the $g$ and $F$ in Theorem 1.1 have the form
$$\displaystyle
g(x,\tau)=\int_0^T e^{-\tau t}\varphi(t)dt\,\frac{\partial v_{\tau}}{\partial\nu}(x),\,\,x\in\partial\Omega,\,\,\tau>0
\tag {1.20}
$$
and
$$\displaystyle
F(x,\tau)=u_f(x,T).
\tag {1.21}
$$
Thus (1.9) is satisfied with the same $\mu$ as (1.17);
(1.18) is satisfied with $\mu'=-1$.
We have to know also an estimation of $C_1$ in (1.4) from above.
From \cite{DL} it follows that
$$\displaystyle
\Vert u(\,\cdot\,,T)\Vert_{L^2(\Omega)}
=O(\Vert f\Vert_{L^2(0,T;H^{-1/2}(\partial\Omega))})
\tag {1.22}
$$
and thus one can choose
$$\displaystyle
C_1=ch_{\Omega}(\omega).
$$
We have also
$$\displaystyle
C_2=ch_D(\omega)<ch_{\Omega}(\omega)=C_3.
$$
Since $C_1+C_3-2C_2=2c(h_{\Omega}(\omega)-h_D(\omega))$, (1.11) becomes (1.18).

{\bf\noindent Remark 1.5.} Note that the choice (1.13) of $z$
satisfying (1.14) goes back to \cite{Ik4} in which an application
of the enclosure method for inverse source problems for the heat
equations are given.  The point of the choice is: the growing
orders of $\vert z\vert$ and $z\cdot z$ as
$\tau\longrightarrow\infty$ are same. See also \cite{IL} for an
application to the so-called inverse heat conduction problem.  In
one-space dimensional case, in \cite{ION} instead of (1.13) the $z$
having the form
$$\displaystyle
z=c\tau\left(1+i\,\sqrt{1-\frac{1}{c^2\tau}}\,\right)
$$
has been used.  In this case $\text{Re}\,(z\cdot z)=\tau$.
Formula (1.19) can be considered as an extension of a result in one-space dimensional case \cite{ION}.
It means that if one could always {\it control} the initial temperature in the process
of all possible measurements at the boundary to be zero, then one can extract the {\it convex hull} of
unknown inclusions.

Since $f$ is {\it complex-valued}, $u_f$ on $\partial\Omega\times\,]0,\,T[$ can not be directly measured
and should be computed from {\it real} data via the formula
$$\displaystyle
u_f=u_{\text{Re}\,f}+iu_{\text{Im}\,f}.
\tag {1.23}
$$
This is a consequence of the zero initial data. Thus the zero
initial data is essential for this procedure.  Note also that
since both $\text{Re}\,f$ and $\text{Im}\,f$ are highly
oscillatory as $\tau\longrightarrow\infty$ with respect to the space
variables, it will be difficult to prescribe those fluxes on the
boundary directly.  Instead one has to make use of the principle
of superposition to {\it compute} the right-hand side of (1.23) on
$\partial\Omega$ from {\it experimental data} which are generated
by finite numbers of independent simpler input fluxes on
$\partial\Omega$. Needles to say, for this procedure the zero
initial data are also essential.

\subsection{The case when $\gamma_0$ is not necessary constant.}

It is possible to to extend Corollary 1.1 to the case
when $\gamma_0$ is {\it not} necessary a constant.
For simplicity of description assume that $\gamma_0-1\in C_0^{\infty}(\Bbb R^n)$.
We construct a special solution of the equation of (1.5)
which has the form
$$\displaystyle
v(x)\sim\frac{e^{x\cdot z}}{\sqrt{\gamma_0}}
$$
as $\tau\longrightarrow\infty$, where $z$ is given by (1.13).

Following \cite{SU}, we make use of the change of the dependent variable formula
(the Liouville transform):
$$\displaystyle
\frac{1}{\sqrt{\gamma_0}}\nabla\cdot\gamma_0\nabla\left(\frac{1}{\sqrt{\gamma_0}}\,\cdot\,\right)=\triangle-V,
\tag {1.24}
$$
where
$$\displaystyle
V=\frac{\triangle\sqrt{\gamma_0}}{\sqrt{\gamma_0}}.
\tag {1.25}
$$

We find the special solution of (1.5) having the form
$$\displaystyle
v=\frac{e^{x\cdot z}}{\sqrt{\gamma_0}}(1+\epsilon_z),
$$
where $\epsilon_z$ is a new unknown function.
It follows from (1.24) and (1.14) that
the equation for $\epsilon=\epsilon_z$ becomes
$$\displaystyle
\left\{\triangle+2z\cdot\nabla-\tau\left(\frac{1}{\gamma_0}-1\right)-V\right\}\epsilon=\tau\left(\frac{1}{\gamma_0}-1\right)+V.
\tag {1.26}
$$

Thus the problem is to construct a solution of equation (1.26)
such that $\epsilon_z\approx 0$ in $\Omega$ as
$\tau\longrightarrow\infty$. In general this is not an easy task
because of the growing factor $\tau$ on the zeroth-order term in
(1.26). However, we found that:  if the virtual slowness $c$ in
(1.13) is {\it sufficiently large and fixed}, then one can
construct such a solution for all large $\tau>>1$ with an arbitrary small
$\epsilon_z$ by using a
combination of the Fourier transform and perturbation methods in
\cite{SU} for the construction of the so-called complex
geometrical optics solutions of the equation
$\nabla\cdot\gamma_0\nabla v=0$.  Its precise description is the following second result.

\proclaim{\noindent Theorem 1.2.}
Let $-1<\delta<0$ and $a,b\in C^{\infty}_0(\Bbb R^n)$.  Given $\eta>0$
there exist positive constants $C_j=C_j(a,b,\Omega,\delta,\eta)$, $j=1,2$ such that:
if $c\ge C_1$ and $\tau\ge C_2$, then $c^2\tau>1$ and there exists a unique $\epsilon_z\in L^2_{\delta}(\Bbb R^n)$
with $z$ given by (1.13) such that
$$\displaystyle
(\triangle+2z\cdot\nabla-\tau a-b)\epsilon_z=\tau a+b\,\,\text{in}\,\Bbb R^n.
\tag {1.27}
$$
Moreover, $\epsilon_z\vert_{\Omega}$ can be identified with a function in $C^1(\overline\Omega)$ and
$$\displaystyle
\Vert\epsilon_z\Vert_{L^{\infty}(\Omega)}+
\Vert\nabla\epsilon_z\Vert_{L^{\infty}(\Omega)}
\le\eta.
\tag {1.28}
$$
\endproclaim

This theorem indicates the important role of
the virtual slowness $c$ when $\gamma_0$ is not necessary
constant. It is not an accessary!  To the best knowledge of the
author this idea, that is, choosing a {\it large} $c$ and fix, never
been pointed out.

Having this theorem, we obtain a result which corresponds to Corollary 1.1.

Let $0<\eta<<1$ and fix a $c\ge C_1$ in Theorem 1.2.
Let $a=(1/\gamma_0-1)$ and $b$ is given by (1.25).
Let $\epsilon_z$ be the solution of (1.27) constructed in Theorem 1.2.
Define
$$\displaystyle
v_{\tau}(x)=\frac{e^{x\cdot z}}{\sqrt{\gamma_0(x)}}(1+\epsilon_z(x)),\,\,x\in\Omega,\,\,\tau\ge C_2.
\tag {1.29}
$$
The function $e^{-\tau t}v_{\tau}(x)$ satisfies the backward parabolic equation
$$\displaystyle
(\nabla\cdot\gamma_0\nabla+\partial_t)(e^{-\tau t}v_{\tau}(x))=0
$$
and its absolute value has the form
$$\displaystyle
\frac{e^{-\tau(t-cx\cdot\omega)}}{\sqrt{\gamma_0(x)}}\vert 1+\epsilon_z(x)\vert
\approx
e^{-\tau(t-cx\cdot\omega)}.
$$
This again supports the name {\it virtual slowness} of $c$.

It is easy to see that the family $(v_{\tau})_{\tau\ge C_2}$ satisfies (1.5), (1.6) and (1.7) with $C_2=ch_D(\omega)$
, (1.8) with $C_3=ch_{\Omega}(\omega)$ and (1.4) with $C_1=ch_{\Omega}(\omega)$.
Thus applying Theorem 1.1 to this case, we obtain the following corollary.

\proclaim{\noindent Corollary 1.2.} Assume that $\gamma_0-1\in
C_0^{\infty}(\Bbb R^n)$.  Fix the virtual slowness as $c=C_1$,
where $C_1$ is just the same as Theorem 1.2. Let $f$ be the
function of $(x,t)\in\partial\Omega\times]0,\,T[$ having a
parameter $\tau>0$ defined by the equation
$$\displaystyle
f(x,t)=\frac{\partial v}{\partial\nu}(x)\varphi(t),
$$
where $v=v_{\tau}$ is given by (1.29) and a real-valued function
$\varphi\in L^2(0,\,T)$
satisfying the condition (1.17) for a $\mu\in\Bbb R$.
Let $u_f=u_f(x,t)$ be the weak solution of (1.1) for $f=f(x,t;\tau)$.
If $T$ satisfies
$$\displaystyle
T>2c(h_{\Omega}(\omega)-h_D(\omega)),
\tag {1.30}
$$
then
$$\displaystyle
\lim_{\tau\longrightarrow\infty}
\frac{1}{2\tau}
\log\left\vert
\int_{\partial\Omega}
\int_0^Te^{-\tau t}\left(-\overline{v_{\tau}(x)}f(x,t;\tau)+u_f(x,t)\gamma_0\frac{\partial\overline{v_{\tau}}}{\partial\nu}(x)\right)dtdS\right\vert
=h_D(\omega).
$$
\endproclaim

In Corollary 1.2 $c=C_1$ and thus (1.30) should be considered as a {\it restriction} on the length of
the time for data collection.  A sufficient condition to ensure (1.30) is
$$\displaystyle
T\ge 2C_1(h_{\Omega}(\omega)+h_{\Omega}(-\omega)).
$$
Note also that the center of coordinates in Theorem 1.2 and Corollary 1.2 is free from $D$ and $\Omega$.

\subsection{Real VS Complex}

Replace the conditions (1.4), (1.6), (1.7) and (1.8) with
the following ones, respectively:
$$\displaystyle
\Vert F(\,\cdot\,,\tau)\Vert_{L^2(\Omega)}
=O(e^{C_1\sqrt{\tau}}\tau^{\mu_1});
\tag {1.31}
$$
$$\displaystyle
\Vert\nabla v_{\tau}\Vert_{L^2(D)}=O(e^{C_2\sqrt{\tau}}\tau^{\mu_2});
\tag {1.32}
$$
$$\displaystyle
\Vert\nabla v_{\tau}\Vert_{L^2(D)}\ge C'' e^{C_2\sqrt{\tau}}\tau^{\mu_3};
\tag {1.33}
$$
$$
\displaystyle
\Vert v_{\tau}\Vert_{H^1(\Omega)}
=O(e^{C_3\sqrt{\tau}}\tau^{\mu_4});
\tag {1.34}
$$
Then instead of (1.12) we have, for any fixed $T>0$ {\it without} (1.11)
and exactly same $g$ satisfying (1.9) and (1.10) for some $\mu, \mu'\in\Bbb R$
$$\displaystyle
\lim_{\tau\longrightarrow\infty}
\frac{1}{2\sqrt{\tau}}
\log
\left\vert\int_{\partial\Omega}
\left(g\overline v_{\tau}-w\gamma_0\frac{\partial\overline v_{\tau}}{\partial\nu}\right)dS\right\vert
=C_2.
\tag {1.35}
$$
Since the proof is simpler than that of Theorem 1.1,
we omit its description as a theorem.  Explicit examples of $v$ satisfying (1.32), (1.33) and (1.34)
in the case when $\gamma_0\equiv 1$ and $n=3$ are the following:
$$\begin{array}{c}
\displaystyle
v(x;\tau,\omega)=e^{\sqrt{\tau} x\cdot\omega},\,x\in\Bbb R^3,\,\,\omega\in S^2,\\
\\
\displaystyle
v(x;\tau,p)=\frac{e^{-\sqrt{\tau}\vert x-p\vert}}{\vert x-p\vert},\,x\in\Bbb R^3\setminus\{p\},\,
\,p\in\Bbb R^3\setminus\overline\Omega,\\
\\
\displaystyle
v(x;\tau,y)=\frac{e^{\sqrt{\tau}\vert x-y\vert}-e^{-\sqrt{\tau}\vert x-y\vert}}{\vert x-y\vert},\,x\in\Bbb R^3\setminus\{y\},\,\,
v(y;\tau,y)=2\tau,\,y\in\Bbb R^3.
\end{array}
$$
These are all real-valued functions and not oscillatory as
$\tau\longrightarrow\infty$. We think that this non oscillatory
character is an advantage in computing the left-hand side of
(1.12).

If $g$ and $F$ are coming from (1.20) and (1.21) for $f$ given by (1.16) for a $\varphi$ satisfying
(1.17) for a $\mu\in\Bbb R$, then one can choose
$$\displaystyle
C_2=\left\{
\begin{array}{lr}
\displaystyle h_{D}(\omega), & \quad\text{if $v=v(x;\tau,\omega)$,}\\
\\
\displaystyle -d_{D}(p), & \quad\text{if $v=v(x;\tau,p)$,}\\
\\
\displaystyle R_{D}(y), & \quad\text{if $v=v(x;\tau,y)$,}
\end{array}
\right.
$$
where $d_*(p)$ and $R_{*}(y)$ denote the distance of $p$ to $*$ and minimum radius
of the open ball that contains $*$ and centered at $y$.
See Lemma 3.1 in \cite{IKII} and Proposition 3.2 in \cite{IKI} for these facts.

By virtue of (1.22) one can choose
$$\displaystyle
C_1=\left\{
\begin{array}{lr}
\displaystyle h_{\Omega}(\omega), & \quad\text{if $v=v(x;\tau,\omega)$,}\\
\\
\displaystyle -d_{\Omega}(p), & \quad\text{if $v=v(x;\tau,p)$,}\\
\\
\displaystyle R_{\Omega}(y), & \quad\text{if $v=v(x;\tau,y)$.}
\end{array}
\right.
$$
Then, formula (1.35) reproduces Theorems 1.2-1.4 in \cite{IKII}.

However, when $\gamma_0$ is not necessary constant, to construct a suitable $v$ one has to solve
the {\it eikonal equation}
$$\displaystyle
\gamma_0\nabla v\cdot\nabla v=1\,\,\text{in}\,\Omega
$$
and corresponding transport equations. However, this is not a
simple matter in general because of the complicated behaviour of
the characteristic curve $x=x(t)$ under suitable initial
conditions on $(x(t),\xi(t))$:
$$\begin{array}{c}
\displaystyle
\frac{dx}{dt}=2\xi,\\
\\
\displaystyle
\frac{d\xi}{dt}=\frac{2}{\sqrt{\gamma_0(x)}}\nabla\left(\frac{1}{\sqrt{\gamma_0(x)}}\right).
\end{array}
$$

Thus Corollary 1.2 suggests that when the background body is isotropic, however, not necessary
homogeneous, the use of complex geometrical optics solutions are {\it better} than that of geometrical
optics solutions.

\section{Proof of Theorem 1.1.}

In what follows for simplicity we wright $v_{\tau}=v$ and $\gamma_0I_n=\gamma_0$.

Let $R_{\gamma}(\tau)$ and $R_{\gamma_0}(\tau)$ denote the Neumann-to-Dirichlet maps on $\partial\Omega$
for the operators $\nabla\cdot\gamma\nabla-\tau$ and $\nabla\cdot\gamma_0\nabla-\tau$, respectively.
We have
$$\displaystyle
R_{\gamma_0}(\tau)\left(\gamma_0\frac{\partial\overline{v}}{\partial\nu}\vert_{\partial\Omega}\right)=\overline{v}\vert_{\partial\Omega},\,\,
R_{\gamma}(\tau)g=p\vert_{\partial\Omega},
$$
where $p$ solves
$$\begin{array}{c}
\displaystyle
(\nabla\cdot\gamma\nabla-\tau)p=0\,\,\text{in}\,\Omega,\\
\\
\displaystyle
\gamma\nabla p\cdot\nu=g(x,\tau)\,\,\text{on}\,\partial\Omega.
\end{array}
\tag {2.1}
$$
Our starting point is the following identity which is an easy consequence of
equations (1.3), (1.5) and integration by parts:
$$\begin{array}{c}
\displaystyle
\int_{\partial\Omega}
\left(g\overline{v}-w\gamma_0\frac{\partial\overline{v}}{\partial\nu}\right)dS
=\int_{\partial\Omega}g(R_{\gamma_0}(\tau)-R_{\gamma}(\tau))
\left(\gamma_0\frac{\partial\overline{v}}{\partial\nu}\vert_{\partial\Omega}\right)dS\\
\\
\displaystyle
+\int_{\Omega}(\gamma-\gamma_0)\nabla\overline{v}\cdot\nabla(w-p)dx
+e^{-\tau T}
\int_{\Omega}F(x,\tau)\overline{v(x)}dx.
\end{array}
\tag {2.2}
$$

Define $\epsilon=w-p$.  It follows from (1.3) and (2.1) that $\epsilon$ solves
$$\begin{array}{c}
\displaystyle
(\nabla\cdot\gamma\nabla-\tau)\epsilon=e^{-\tau T}F(x,\tau)\,\,\text{in}\,\Omega\\
\\
\displaystyle
\gamma\nabla\epsilon\cdot\nu=0\,\,\text{on}\,\partial\Omega.
\end{array}
$$
Since $\tau>0$, it is easy to see that
$$\displaystyle
\Vert\nabla\epsilon\Vert_{L^2(\Omega)}
\le C e^{-\tau T}\tau^{-1/2}
\Vert F(\,\cdot\,,\tau)\Vert_{L^2(\Omega)}
$$
and from (1.4) one gets
$$\displaystyle
\Vert\nabla\epsilon\Vert_{L^2(\Omega)}
=O(e^{-\tau(T-C_1)}\tau^{\mu_1-1/2}).
$$
From this together with (1.4), (1.6) and (1.8) we obtain
$$\begin{array}{c}
\displaystyle
\int_{\Omega}(\gamma-\gamma_0)\nabla\overline{v}\cdot\nabla\epsilon dx
+e^{-\tau T}\int_{\Omega}F(x,\tau)\overline{v(x)}dx\\
\\
\displaystyle
=O(e^{C_2\tau}\tau^{\mu_2}e^{-\tau(T-C_1)}\tau^{\mu_1-1/2})
+O(e^{-\tau T}e^{C_3\tau}\tau^{\mu_4}\tau^{\mu_1}e^{C_1\tau})\\
\\
\displaystyle
=O(e^{-\tau(T-C_1-C_2)}\tau^{\mu_1+\mu_2-1/2})
+O(e^{-\tau(T-C_1-C_3)}\tau^{\mu_1+\mu_4}).
\end{array}
\tag {2.3}
$$

A combination of (2.2) and (2.3) gives
$$\begin{array}{c}
\displaystyle
\int_{\partial\Omega}
\left(g\overline{v}-w\gamma_0\frac{\partial\overline{v}}{\partial\nu}\right)dS
=\Psi(\tau)
\int_{\partial\Omega}\gamma_0\frac{\partial v}{\partial\nu}(R_{\gamma_0}(\tau)-R_{\gamma}(\tau))\left(\gamma_0\frac{\partial \overline{v}}{\partial\nu}\vert_{\partial\Omega}\right)dS\\
\\
\displaystyle
+O(e^{-\tau(T-C_1-C_2)}\tau^{\mu_1+\mu_2-1/2})
+O(e^{-\tau(T-C_1-C_3)}\tau^{\mu_1+\mu_4}).
\end{array}
\tag {2.4}
$$

The following type of estimates now are well known and it is a
consequence of Proposition 2.1 in \cite{IKII} which goes back to
\cite{IES} and one of assumptions (A1) and (A2).
See also \cite{KSS} when $\gamma_0=1$ and $h=(k-1)I_n$ with a positive constant $k$.

\proclaim{\noindent Lemma 2.1.}
There exist $C>0$ and $C'>0$ such that for all $v$ satisfying (1.5)
$$\displaystyle
C\Vert\nabla v\Vert_{L^2(D)}^2\le
\left\vert\int_{\partial\Omega}\gamma_0\frac{\partial v}{\partial\nu}(R_{\gamma_0}(\tau)-R_{\gamma}(\tau))
\left(\gamma_0\frac{\partial\overline{v}}{\partial\nu}\vert_{\partial\Omega}\right)dS\right\vert
\le C'\Vert\nabla v\Vert_{L^2(D)}^2.
\tag {2.5}
$$

\endproclaim

Note also that:

$\bullet$  (1.9) is equivalent to the statement: there exists a $C_0>0$ such that,
as $\tau\longrightarrow\infty$
$$\displaystyle
C_0\le\tau^{\mu}\vert\Psi(\tau)\vert.
\tag {2.6}
$$

From the right-hand side of (2.5), (1.6), (2.4), (1.10) and (2.3)
we have
$$\begin{array}{c}
\displaystyle
\left\vert\int_{\partial\Omega}
\left(g\overline{v}-w\gamma_0\frac{\partial\overline{v}}{\partial\nu}\right)dS\right\vert\\
\\
\displaystyle
=O(\tau^{\mu'} e^{2C_2\tau}\tau^{2\mu_2})
+O(e^{-\tau(T-C_1-C_2)}\tau^{\mu_1+\mu_2-1/2})
+O(e^{-\tau(T-C_1-C_3)}\tau^{\mu_1+\mu_4})\\
\\
\displaystyle
=O(e^{2C_2\tau}\tau^{2\mu_2+\mu'}
(1+e^{-\tau(T-C_1+C_2)}\tau^{\mu_1-\mu_2-\mu'-1/2}
+e^{-\tau(T-C_1+2C_2-C_3)}\tau^{\mu_1+\mu_4-2\mu_2-\mu'})).
\end{array}
\tag {2.7}
$$

One the other hand, from the left-hand side of (2.5), (2.6) and (1.7)
we obtain
$$\begin{array}{c}
\displaystyle
\left\vert\int_{\partial\Omega}
\left(g\overline{v}-w\gamma_0\frac{\partial\overline{v}}{\partial\nu}\right)dS\right\vert\\
\\
\displaystyle
\ge CC_0(C'')^2\tau^{-\mu}e^{2C_2\tau}\tau^{2\mu_3}
+O(e^{-\tau(T-C_1-C_2)}\tau^{\mu_1+\mu_2-1/2})
+O(e^{-\tau(T-C_1-C_3)}\tau^{\mu_1+\mu_4})\\
\\
\displaystyle
=C'''e^{2C_2\tau}\tau^{2\mu_3-\mu}
(1+ O(e^{-\tau(T-C_1+C_2)}\tau^{\mu_1+\mu_2-2\mu_3+\mu-1})
+O(e^{-\tau(T-C_1+2C_2-C_3)}\tau^{\mu_1+\mu_4-2\mu_3+\mu})),
\end{array}
\tag {2.8}
$$
where $C'''=CC_0(C'')^2>0$.

Now formula (1.12) is a consequence of (1.11), (2.7) and (2.8)
provided
$$\displaystyle
T>\max(C_1-C_2,C_1+C_3-2C_2).
\tag {2.9}
$$
However, by Remark 1.2 we have $(C_1+C_3-2C_2)-(C_1-C_2)=C_3-C_2\ge 0$ and thus
(2.9) is nothing but (1.11).

\noindent
$\Box$

\section{Proof of Theorem 1.2}

\subsection{A special fundamental solution}

Given $F$ we construct a solution of the inhomogeneous modified Helmholtz equation
$$\displaystyle
(-\triangle+\tau)v+F=0\,\,\text{in}\,\Bbb R^n
\tag {3.1}
$$
that has the form
$$\displaystyle
v(x)=e^{x\cdot z}\Psi(x)
\tag {3.2}
$$
where the complex vector $z$ is given by (1.13).

Write
$$\displaystyle
F(x)=e^{x\cdot z}f(x).
\tag {3.3}
$$
Then if $\Psi$ satisfies the equation
$$\displaystyle
-\triangle\Psi-2z\cdot\nabla\Psi+f=0\,\,\text{in}\,\Bbb R^n,
\tag {3.4}
$$
then the $v$ given by (3.2) satisfies (3.1) with $F$ given by (3.3).

We construct a solution of (3.4) by using a special fundamental
solution of (3.4).
Set
$$\displaystyle
Q_z(\xi)=\vert\xi\vert^2-2iz\cdot\xi,\,\,\xi\in\Bbb R^n.
$$
We have
$$\displaystyle
\text{Re}\,Q_z(\xi)=\left\vert\xi+c\tau\sqrt{1-\frac{1}{c^2\tau}}\,\omega^{\perp}\right\vert^2
-c^2\tau^2\left(1-\frac{1}{c^2\tau}\right)
$$
and
$$\displaystyle
\text{Im}\,Q_z(\xi)=-2c\tau\omega\cdot\xi.
$$
Thus the set $Q_{z}^{-1}(0)=\{\xi\in\Bbb R^n\,\vert\,Q_z(\xi)=0\}$ consists of
the circle on the plane $\omega\cdot\xi=0$ centered at $-c\tau\sqrt{1-(c^2\tau)^{-1}}\omega^{\perp}$
with radius $c\tau\sqrt{1-(c^2\tau)^{-1}}$. Moreover $\nabla\text{Re}\,Q_z(\xi)$ and $\nabla\text{Im}\,Q_z(\xi)$
on $Q_z^{-1}(0)$ are linearly independent.  Thus $-1/Q_z(\xi)$ is locally integrable on $\Bbb R^n$
and can be identified with a unique tempered distribution and its inverse Fourier transform
$$\displaystyle
g_z(x)=-\frac{1}{(2\pi)^n}
\int_{\Bbb R^n}\frac{e^{ix\cdot\xi}}{Q_z(\xi)}d\xi
\tag {3.5}
$$
is well defined.  This $g=g_z$ satisfies, in the sense of tempered distribution
$$\displaystyle
-\triangle g-2z\cdot\nabla g+\delta(x)=0\,\,\text{in}\,\Bbb R^n
$$
and thus, in the sense of Schwartz distribution
$$\displaystyle
(-\triangle+\tau)(e^{x\cdot z}g_z)+\delta(x)=0\,\,\text{in}\,\Bbb R^n.
$$
The $G_z(x)=e^{x\cdot z}g_z(x)$ should be called the Faddeev-Green function.

From \cite{SU}, one knows that: by virtue of Theorem 7.1.27 of \cite{H} given
$f\in L_{\delta+1}(\Bbb R^n)$ with $-<1<\delta<0$ the solution of
equation (3.4) is unique in $L^2_{\delta}(\Bbb R^n)$; if $f$ is a
rapidly decreasing function, then the unique solution of (3.4) is
given by
$$\displaystyle
\Psi=g_z*f;
\tag {3.6}
$$
its first and second derivatives belong to $L^2_{\delta}(\Bbb R^n)$
and satisfy
$$\displaystyle
\Vert D^{\alpha}\Psi\Vert_{L^2_{\delta}(\Bbb R^n)}
\le C_{\alpha}(z,\delta)\Vert f\Vert_{L^2_{1+\delta}(\Bbb R^n)},\,\,\vert\alpha\vert\le 2.
\tag {3.7}
$$
Thus the operator :$f\longmapsto \Psi$ has the unique continuous
extension as a bounded linear operator of $L^2_{\delta+1}(\Bbb
R^n)$ to $L^2_{\delta}(\Bbb R^n)$.  We still denote the operator
by the same symbol $g_z*f$ which yields the unique solution of
equation (3.4) in $L^2_{\delta}(\Bbb R^n)$ for {\it all} $f\in
L^2_{\delta+1}(\Bbb R^n)$.

However, for our purpose, one has to clarify the behaviour of
$C_{\alpha}(z,\delta)$ in (3.7) as $\tau\longrightarrow\infty$.
This is not the well known case $z\cdot z=-k^2$ with
a $k\ge 0$ which appeared in inverse scattering/boundary
value problems (cf.\cite{SU, N1, No})
since we have (1.14) and thus $z\cdot z\longrightarrow\infty$
as $\tau\longrightarrow\infty$.
In the next subsection we study
more about $C_{\alpha}(z,\delta)$ by carefully checking a proof of
(3.7).

\subsection{Asymptotic behaviour}

Define
$$\displaystyle
\lambda(c,\tau)=\sqrt{1-\frac{1}{c^2\tau}}.
$$

The aim of this subsection is to give a proof of the following estimates.

\proclaim{\noindent Proposition 3.1.}
Let $R>0$ and $-1<\delta<0$.
We have
$$\displaystyle
\Vert D^{\alpha}g_{z}*f\Vert_{L^2_{\delta}(\Bbb R^n)}
\le
(c\tau\lambda)^{\vert\alpha\vert-1}C_{\delta,R}
\Vert f\Vert_{L^2_{\delta+1}(\Bbb R^n)},\,\vert\alpha\vert\le 2,\,c\tau\lambda\ge R.
\tag {3.8}
$$

\endproclaim

A change of variables yields
$$\displaystyle
g_z(x)
=(c\tau\lambda(c,\tau))^{n-2}h_{\lambda}(y;\omega,\omega^{\perp})\vert_{y=c\tau\lambda(c,\tau)x,\,\lambda=\lambda(c,\tau)}
\tag {3.9}
$$
where
$$\displaystyle
h_{\lambda}(y;\omega,\omega^{\perp})
=-\frac{1}{(2\pi)^n}\int_{\Bbb R^n}\frac{e^{iy\cdot\xi}d\xi}
{\vert\xi+\omega^{\perp}\vert^2-1-2i\lambda^{-1}\omega\cdot\xi},\,0<\lambda\le 1.
\tag {3.10}
$$

Given an arbitrary rapidly decreasing function $f$ we give an estimate
of $h_{\lambda}*f$.
Using a rotation in $\xi$-space, it suffices to consider the case when $\omega=(0,1,0,\cdots,0)^T$
and $\omega^{\perp}=(1,0,\cdots,0)^T$.

Define
$$\displaystyle
\Sigma=\{\xi\in\Bbb R^n\,\vert\,\vert\xi+\omega^{\perp}\vert^2=1,\,\,\omega\cdot\xi=0\}.
$$
$\Sigma$ is the set of all real zero points of the complex-valued polynomial
$$\displaystyle
h(\xi;\lambda)=\vert\xi+\omega^{\perp}\vert^2-1-2i\lambda^{-1}\omega\cdot\xi,
$$
however, $\Sigma$
itself is {\it independent} of $\lambda$.

\proclaim{\noindent Lemma 3.1.}  Given $\epsilon>0$ there exists a positive constant $C_{\epsilon}$ independent
of $\lambda$ such that for all $\xi\in\Bbb R^n$ with $\text{dist}_{\Sigma}(\xi)\ge\epsilon$,
we have
$$\displaystyle
\vert h(\xi;\lambda)\vert\ge C_{\epsilon}\vert\xi\vert^2.
$$

\endproclaim

{\it\noindent Proof.}
Since $0<\lambda\le 1$, we have
$$\displaystyle
\vert h(\xi;\lambda)\vert\ge\vert h(\xi;1)\vert.
$$
Let $\vert\xi\vert\ge 2$.  We have
$$\begin{array}{c}
\displaystyle
\vert h(\xi;1)\vert^2
=(\vert\xi\vert^2+2\omega^{\perp}\cdot\xi)^2+4(\omega\cdot\xi)^2\\
\\
\displaystyle
=\vert\xi\vert^4+4(\omega^{\perp}\cdot\xi)\vert\xi\vert^2+4(\omega^{\perp}\cdot\xi)^2
+4(\omega\cdot\xi)^2\\
\\
\displaystyle
\ge\vert\xi\vert^4-4\vert\xi\vert^3
=\vert\xi\vert^4\left(1-\frac{1}{\vert\xi\vert}\right)
\ge\frac{\vert\xi\vert^4}{2}.
\end{array}
$$
Since $\xi\not=0$ and $\vert h(\xi;1)\vert>0$ for $\xi$ with $\text{dist}\,_{\Sigma}(\xi)\ge\epsilon$, the function
$\vert\xi\vert^{-2}\vert h(\xi;1)\vert$ is continuous and positive on the nonempty compact set
$K_{\epsilon}=\{\xi\in\Bbb R^n\,\vert\, \vert\xi\vert\le 2\,\text{and}\,\text{dist}\,_{\Sigma}(\xi)\ge\epsilon\}$.
Thus $m_{\epsilon}\equiv\inf_{\xi\in K_{\epsilon}}\vert\xi\vert^{-2}\vert h(\xi;1)\vert>0$ and choosing
$$\displaystyle
C_{\epsilon}=\min\left(\frac{1}{\sqrt{2}},m_{\epsilon}\right),
$$
we obtain the desired estimate.

\noindent
$\Box$

For the treatment of $h(\xi;\lambda)$ in a neighbourhood of $\Sigma$,
we start with the following fundamental fact:
$$\displaystyle
-\frac{1}{2\pi i(x_1+ix_2)}=\frac{1}{(2\pi)^2}\int\frac{e^{ix\cdot\eta}}{\eta_1+i\eta_2}d\eta.
$$
This yields
$$\displaystyle
-\frac{1}{2\pi i(\lambda^{-1}x_1+ix_2)}
=\frac{1}{(2\pi)^2}\int\frac{e^{ix\cdot\eta}}{\eta_1+i\lambda^{-1}\eta_2}d\eta.
\tag {3.11}
$$

\noindent

\proclaim{\noindent Lemma 3.2.} Let $-1<\delta<0$ and $0<\lambda\le 1$.  Then
$$\displaystyle
Z_{\lambda}f=\check{\left(\frac{\hat{f}}{\eta_1+i\lambda^{-1}\eta_2}\right)}
$$
defines a bounded map from $L^2_{\delta+1}(\Bbb R^n)$ to $L^2_{\delta}(\Bbb R^n)$
and its operator norm is bounded from above with respect to $\lambda$.
\endproclaim

{\bf\noindent Proof.}
From (3.11) we have
$$\displaystyle
Z_{\lambda}f=-\frac{1}{2\pi i}\left\{\frac{1}{\lambda^{-1}x_1+ix_2}\right\}\ast f,
$$
where $*$ denotes convolution with respect to variables
$(x_1,x_2)$. Let $g$ be an arbitrary rapidly decreasing function.
Since $0<\lambda\le 1$ we have
$$\displaystyle
\vert\lambda^{-1}x_1+ix_2\vert\ge\vert x\vert.
$$
This gives
$$\begin{array}{c}
\displaystyle
(2\pi)^2\vert <Z_{\lambda}f,g>\vert^2
\le
\left(\int\int\left\vert\frac{g(x)f(y)}{\lambda^{-1}(x_1-y_1)+i(x_2-y_2)}\right\vert dxdy\right)^2\\
\\
\displaystyle
\le\left(\int\int\left\vert\frac{g(x)f(y)}{(x_1-y_1)+i(x_2-y_2)}\right\vert dxdy\right)^2.
\end{array}
$$
Thus everything is reduced to the case when $\lambda=1$
and it is nothing but Lemma 3.1 in \cite{SU}.

\noindent
$\Box$

Having Lemmas 3.1-2 and using a rotation in $\xi$-space, localization and a change of variable which
are exactly same as Sylvester-Uhlmann's argument \cite{SU} for the case $z\cdot z=0$,
we obtain
$$\displaystyle
\Vert D^{\alpha}h_{\lambda}*f\Vert_{L^2_{\delta}(\Bbb R^n)}
\le C_{\delta}\Vert f\Vert_{L^2_{\delta+1}(\Bbb R^n)},\,\vert\alpha\vert\le 2.
\tag {3.12}
$$
The constant $C_{\delta}$ is independent of $\lambda$.

To deduce (3.8) from (3.12) we employ a scaling argument \cite{Ik0}.
It follows from (3.9) that
$$\displaystyle
g_z*f
=(c\tau\lambda)^{-2}\left(h_{\lambda}* f_{(c\tau\lambda)^{-1}}\right)_{c\tau\lambda},\,\,\lambda=\lambda(c,\tau),
\tag {3.13}
$$
where
$$\displaystyle
f_{\eta}(x)=f(\eta\,x),\,\,\eta>0.
$$
Let $s>0$.  We have
$$\displaystyle
\Vert f_{\eta^{-1}}\Vert_{L^2_s(\Bbb R^n)}\le \eta^{s+n/2}\max\,(1,R^{-s})\Vert f\Vert_{L^2_s(\Bbb R^n)},\,0<R\le \eta.
$$
Now from this and (3.13) we have
$$\begin{array}{c}
\displaystyle
\vert<g_z\ast f,g>\vert=(c\tau\lambda)^{-2}\vert<(h_{\lambda}*f_{(c\tau\lambda)^{-1}})_{c\tau\lambda}, g>\vert\\
\\
\displaystyle
=(c\tau\lambda)^{-(2+n)}\vert<h_{\lambda}* f_{(c\tau\lambda)^{-1}},g_{(c\tau\lambda)^{-1}}>\vert\\
\\
\displaystyle
\le (c\tau\lambda)^{-(2+n)}\Vert h_{\lambda}*f_{(c\tau\lambda)^{-1}}\Vert_{L^2_{\delta}(\Bbb R^n)}\Vert g_{(c\tau\lambda)^{-1}}\Vert_{L^2_{-\delta}(\Bbb R^n)}
\\
\\
\displaystyle
\le (c\tau\lambda)^{-(2+n)}C_{\delta}\Vert f_{(c\tau\lambda)^{-1}}\Vert_{L^2_{\delta+1}(\Bbb R^n)}\Vert g_{(c\tau\lambda)^{-1}}\Vert_{L^2_{-\delta}(\Bbb R^n)}\\
\\
\displaystyle
\le (c\tau\lambda)^{-(2+n)}C_{\delta}(c\tau\lambda)^{\delta+1+n/2}\max(1, R^{-(\delta+1)})
\Vert f\Vert_{L^2_{\delta+1}(\Bbb R^n)}
(c\tau\lambda)^{-\delta+n/2}\max(1,R^{\delta})\Vert g\Vert_{L^2_{-\delta}(\Bbb R^n)}\\
\\
\displaystyle
=(c\tau\lambda)^{-1}C(\delta,R)\Vert f\Vert_{L^2_{\delta+1}(\Bbb R^n)}\Vert g\Vert_{L^2_{-\delta}(\Bbb R^n)}.
\end{array}
$$
This together with the same argument for $D^{\alpha}g_z*f$ yields (3.8).

\subsection{Uniqueness and Construction of $\epsilon_z$}

Write (1.24) as
$$\displaystyle
(-\triangle-2z\cdot\nabla)\epsilon_z=(\tau a+b)\epsilon_z-(\tau a+b)\,\,\text{in}\,\Bbb R^n.
$$
Since the solution of (3.4) is unique in $L^2_{\delta}(\Bbb R^n)$ and has the form (3.6) for $f\in L^2_{\delta+1}(\Bbb R^n)$,
we have
$$\displaystyle
\epsilon_z=\tau g_z*(a\epsilon_z)-\tau g_z*a+g_z*(b\epsilon_z)-g_z*b.
\tag {3.14}
$$
Define
$$\displaystyle
K_z(a):L^2_{\delta}(\Bbb R^n)\ni h\longmapsto \tau g_z*(ah)\in L^2_{\delta}(\Bbb R^n)
$$
and
$$\displaystyle
L_z(b):L^2_{\delta}(\Bbb R^n)\ni h\longmapsto g_z*(bh)\in L^2_{\delta}(\Bbb R^n).
$$
It follows from (3.7) that both $K_z$ and $L_z$ define bounded linear operators
in $L^2_{\delta}(\Bbb R^n)$.  Rewrite (3.14) as
$$\displaystyle
(I-K_z(a)-L_z(b))\epsilon_z=-\tau g_z*a-g_z*b.
$$
It follows from (3.8) for $\vert\alpha\vert=0$ that
$$\displaystyle
\Vert K_z(a)h\Vert_{L^2_{\delta}(\Bbb R^n)}\le
(c\lambda)^{-1}C_{\delta,R}
\Vert ah\Vert_{L^2_{\delta+1}(\Bbb R^n)}
$$
provided $c\lambda\ge R>0$.
Since $\Vert ah\Vert_{L^2_{\delta+1}(\Bbb R^n)}\le\Vert<x>a\Vert_{L^{\infty}(\Bbb R^n)}
\Vert h\Vert_{L^2_{\delta}(\Bbb R^n)}$, one gets
$$\displaystyle
\Vert K_z(a)\Vert\le
(c\lambda)^{-1}C_{\delta,R}\Vert<x>a\Vert_{L^{\infty}(\Bbb R^n)}
$$
and similarly
$$\displaystyle
\Vert L_z(b)\Vert\le
(c\tau\lambda)^{-1}C_{\delta,R}\Vert<x>b\Vert_{L^{\infty}(\Bbb R^n)}.
$$
From these we obtain
$$\displaystyle
\Vert K_z(a)+L_z(b)\Vert\le(c\lambda)^{-1}C(\delta,R)
(\Vert<x>a\Vert_{L^{\infty}(\Bbb R^n)}+\tau^{-1}\Vert<x>b\Vert_{L^{\infty}(\Bbb R^n)}).
\tag {3.15}
$$

Let $c_1\ge c_2>0$ and $R>0$.
Let $c$ and $\tau$ satisfy
$$\displaystyle
c\ge\sqrt{c_1^2+R^2}
$$
and
$$\displaystyle
\tau\ge\frac{1}{c_2^2}.
$$
We have $c^2\tau>1$ and $c\lambda\ge R$.  Now given $0<\delta<1$, $R>0$ and $\eta>0$
choose a large $c_1$ in such a way that
$$\displaystyle
\eta\sqrt{c_1^2-c_2^2}
\ge C(\delta,R)
(\Vert<x>a\Vert_{L^{\infty}(\Bbb R^n)}+c_2^2\Vert<x>b\Vert_{L^{\infty}(\Bbb R^n)}).
\tag {3.16}
$$
Since $c\lambda\ge\sqrt{c_1^2-c_2^2}$, a combination of (3.15) and (3.16) gives
$$\displaystyle
\Vert K_z(a)+L_z(b)\Vert\le\eta.
\tag {3.17}
$$
Choosing a $\eta<1$ and
$$\displaystyle
C_1=\sqrt{c_1^2+R^2},\,\,C_2=1/c_2^2
$$
for $c_1$ satisfying (3.16), specially chosen $c_2\le c_1$ and $R$,
we obtain the uniqueness and existence of the solution of (3.12) and thus those of (1.27), too.

$\epsilon_z\in L^2_{\delta}(\Bbb R^n)$ takes the form
$$\displaystyle
\epsilon_z
=-\sum_{n=0}^{\infty}(K_z(a)+L_z(b))^n(\tau g_z*a+g_z*b).
$$
It follows from (3.17) and (3.8) for $\vert\alpha\vert=0$ that
$$\begin{array}{c}
\displaystyle
\Vert\epsilon_z\Vert_{L^2_{\delta}(\Bbb R^n)}
\le\sum_{n=0}^{\infty}\eta^n
(c\lambda)^{-1}C(\delta,R)
(\Vert a\Vert_{L^{2}_{\delta+1}(\Bbb R^n)}+c_2^2\Vert b\Vert_{L^{2}_{\delta+1}(\Bbb R^n)})\\
\\
\displaystyle
=\frac{1}{1-\eta}(c\lambda)^{-1}C(\delta,R)
(\Vert a\Vert_{L^{2}_{\delta+1}(\Bbb R^n)}+c_2^2\Vert b\Vert_{L^{2}_{\delta+1}(\Bbb R^n)})
=O((c\lambda)^{-1}).
\end{array}
\tag {3.18}
$$
Differentiating both sides of (3.14) in the sense of distribution, we have
$$\displaystyle
(I-K_z(a)-K_z(b))\partial_j\epsilon_z
=-\tau g_z*\partial_{j}a-K_z(\partial_{j}a)\epsilon_z
-g_z*\partial_{j}b-L_z(\partial_{j}b)\epsilon_z.
$$
A similar argument for the derivation of (3.18) yields
$$\displaystyle
\Vert\partial_j\epsilon_z\Vert_{L^2_{\delta}(\Bbb R^n)}
=O((c\lambda)^{-1}).
$$
Continue this procedure for higher order derivatives of
$\epsilon_z$ in $L^2_{\delta}(\Bbb R^n)$ and apply the Sobolev
imbedding theorem to the resulted estimates. Then one concludes
that $\epsilon_z\vert_{\Omega}$ can be identified with a function
in $C^1(\overline\Omega)$ and
$$\displaystyle
\Vert\epsilon_z\Vert_{L^{\infty}(\Omega)}
+\Vert\nabla\epsilon_z\Vert_{L^{\infty}(\Omega)}
=O((c\lambda)^{-1}).
$$
Thus choosing again a large $c_1$, we obtain (1.28).
This completes the proof of Theorem 1.2.

$$\quad$$

\centerline{{\bf Acknowledgement}}

This research was partially supported by Grant-in-Aid for
Scientific Research (C)(No.  21540162) of Japan  Society for
the Promotion of Science.

$$\quad$$

\vskip1cm
\noindent
e-mail address

ikehata@math.sci.gunma-u.ac.jp

\end{document}